# RENORMALIZATION AND BLOW UPS FOR THE NONLINEAR SCHRÖDINGER EQUATION

D. GAIDASHEV

ABSTRACT. Existence of finite-time blow ups in the classical one-dimensional nonlinear Schrödinger equation (NLS)

$$(1) \qquad i\partial_t u + u_{xx} + |u|^{2r} u = 0, \qquad u(x,0) = u_0(x)$$

has been one of the central problems in the studies of the singularity formation in the PDE's.

We revisit this problem using a novel approach based on the ideas borrowed from Dynamical Systems.

To that end, we reformulate the initial value problem for (1), with $r \in \mathbb{N}$, $r \geq 1$, as a fixed point problem for a certain *renormalization* operator, and use the ideas of apriori bounds to prove existence of a renormalization fixed point. Existence of such fixed points leads to existence of self-similar solutions of the form

$$u(x,t) = (T-t)^{-\frac{1}{2r}} U\left((T-t)^{-\frac{1}{2}} x\right),$$

whose $L^{2r+2}$-norms are bounded up-to a finite time $T$ and whose energy blows up at $T$.

## Contents

## 1. Introduction

In this paper we study existence of self-similar blow up solutions in the 1D non-linear Schrödinger equation with a generalizied power non-linearity

$$(2) \qquad i\partial_t u + \Delta u + \lambda |u|^{2r} u = 0, \qquad u(x,0) = u_0(x), \qquad x \in \mathbb{R},$$

where $r \in \mathbb{N}$, $r \geq 1$.

---

*Date*: April 9, 2025.

1.1. **Global existence results.** Global existence of solutions of (2) has been the subject of many works. We will now list the most important local and global existence results obtained in the literature.

Let $d$ be the dimension of the underlying real space, $\mathbb{R}^d$. Denote

(3) $$\pi(d) = \begin{cases} \infty, & d = 1, 2, \\ \frac{2}{d-2}, & d \geq 3. \end{cases}$$

- [1] demonstrates that if $1 \leq d \leq 3$, $0 < r < \pi(d)$, $\lambda < 0$, the IVP for (2) with the initial data $u_0 \in H^2(\mathbb{R}^d)$ has a unique global solution $u(t) \in C(\mathbb{R}, H^2(\mathbb{R}^d)) \cap C^1(\mathbb{R}, L^2(\mathbb{R}^d))$.

- It has been shown in [7] that in the case of $d \geq 2$, $0 < r < \pi(d)$, $\lambda < 0$, or $d \geq 2$, $0 < r < 2/d$, $\lambda > 0$, equation (2) has a globally defined weak solution $u(t) \in C(\mathbb{R}, H^1(\mathbb{R}^d))$ whenever $u_0(x) \in H^1(\mathbb{R}^d)$.

- [19] shows that in the case of $r > 0$, $\lambda < 0$, the NLS (2) has at least one globally defined weak solution $u(t) \in L^\infty(\mathbb{R}, H^1(\mathbb{R}^d) \cap L^{2r+2}(\mathbb{R}^d))$ whenever $u_0(x) \in H^1(\mathbb{R}^d) \cap L^{2r+2}(\mathbb{R}^d)$.

- [18] shows that in the case of $r = 2/d$, $\lambda > 0$, the NLS (2) has unique locally defined weak solution $u(t) \in L^\infty((-\infty, T), L^2(\mathbb{R}^d)) \cap L^{2r+2}((-\infty, T) \times \mathbb{R}^d)$ whenever $u_0(x) \in L^2(\mathbb{R}^d)$.

- [20] shows that in the case of $0 < r < 2/d$, $\lambda > 0$, the NLS (2) has unique global solution $u(t) \in C(\mathbb{R}, L^2(\mathbb{R}^d)) \cap L_{loc}^{\frac{4(r+1)}{dr}}(\mathbb{R}, L^{2r+2}(\mathbb{R}^d))$ whenever $u_0(x) \in L^2(\mathbb{R}^d)$.

1.2. **Results regarding singularity formation.** Conversely, a number of works addressed the issue of existence of singular solutions.

- Existence of singular solution of (2) for $r \geq 2/d$ has been demonstrated in [8] and [10]. Specifically, it is proved in [8] that if $u_0(x)$ belongs to a Schwartz class, if the initial energy is nonpositive and if $\text{Im} \int \rho \overline{u_0(x)} \partial_\rho u_0(x) dx > 0$, where $\rho = |x|$, then the energy of this solution blows up in finite time.

- The NLS with $r = 2$ has a special property of conformal invariance, that is, if $u(x, t)$ is a solution, then so is

$$(T-t)^{-\frac{1}{2}} e^{\frac{-i|x|^2}{4(T-t)}} \overline{u(x(T-t)^{-1}, (T-t)^{-1})}$$

for $t < T$. Here $\overline{u}$ denotes the complex conjugation. This allows one to construct a very specific explicit blow up solution as follows.

The unique even solution $Q(x)$ to the problem

(4) $$\begin{cases} u = \Delta u + |u|^4 u \\ u > 0 \text{ in } \mathbb{R} \end{cases}$$

(see [4], [3], [17]), [11]) is usually called the ground state, and plays an important role in the case of $r = 2$; if $\|u_0\|_2 < \|Q\|_2$ then the solution $u(x, t)$ to the IVP does not blow up; if $\|u_0\|_2 \geq \|Q\|_2$, it may blow up in finite time. In particular, any blow up with $\|u_0\|_2 = \|Q\|_2$ has the form

$$u(x,t) = \frac{\omega^{\frac{1}{2}}}{(T-t)^{\frac{1}{2}}} e^{-i\frac{|x|^2 - 4w^2}{4(T-t)} + i\theta} Q\left(\frac{\omega x}{T-t}\right), \tag{5}$$

where $\omega$ and $\theta$ are parameters, see [16].

1.3. **Renormalization in evolution PDEs.** A novel approach to existence of singular solutions to nonlinear PDEs has been proposed by Bakhtin, Dinaburg and Sinai in [2]. In [2] a solution to an IVP with a certain self-similar initial conditions is "shadowed" by a solution of a fixed point problem for an integral nonlinear "renormalization" operator in an appropriate functional space. The self-similarity of the solution is of the type

$$u(x,t) = (T-t)^{-\frac{1}{2r}} U\left((T-t)^{-\frac{1}{2}} x\right), \tag{6}$$

different from (5).

Later this renormalization approach was applied to several hydrodynamics models by Li and Sinai [12–15]. The method of these publications is a development of that of [2], with the addition that the authors derive an exact, and not an approximate, renormalization fixed point equation.

We have used the renormalization point of view of [2] to prove existence of blow up solution in the 1D geostrophic equation in [6], and existence of certain "non-physical" blow up solutions to the generalized Navier-Stokes equation in [5].

In this paper we will use technique similar to those in [6] and [5] to first derive an equation for the Fourier transform of a certain self-similar solution to NLS and demonstrate that this equation reduces to a fixed point problem for a renormalization operator. We then prove existence of a family of exponentially decaying renormalization fixed points, which correspond to real $C^\infty$-solutions of the Picard-Lindelöff equation for (2) which blow up in finite time $T$.

We would like to underline that the self-similarity of the blow ups is different from that of the minimal mass solution $Q$; this self-similarity is of form (6) and is present for all integer $r \geq 1$. The blow up solution itself is not an $L^2$-function, but rather its $L^{2r+2}$-norm, together with the $L^2$-norm of the first spacial derivative, is bounded up to the blow up time.

1.4. **Outline of the proof.** As usual, a solution to a Picard-Lindelöff integral equation associated to a PDE, will be called a *mild* solution. Specifically, we will says that $u(x,t)$ is a mild solution to the IVP (2) if for all $x \in \mathbb{R}$ and $t \in [0, T)$, where $T$ can be finite or $T = \infty$,

$$u(x,t) = S(t)u_0(x) + i\int_0^t S(t-s)|u(x,s)|^{2r} u(x,s) \mathrm{d}s, \tag{7}$$

where $S$ is the Schrödinger group with the infinitesimal generator $i\Delta$.

In Section 3 we will demonstrate that the problem of existence of mild self-similar blow ups is equivalent to a problem of existence of a common fixed point of a family of integral operators $\mathcal{R}_\beta$ acting on functions of one variable in a weighted Lebegue space, with $\beta \in (0,1)$ playing the role of "time": $\beta^2 = (T-t)/T$.

In Section 4 this fixed point problem will be reduced to a construction of relatively compact sets $\mathcal{N}_p$ in a weighted Lebesgue space, invariant by all $\mathcal{R}_\beta$ with $\beta \in (\beta_0, 1)$ for some $\beta_0 \in (0,1)$. The sets $\mathcal{N}_p$ will be constructed as the non-empty intersections of two convex sets $\tilde{\mathcal{E}}_D$ and $\mathcal{M}_\mu$. The first, $\tilde{\mathcal{E}}_D$, is the set of all $\psi$ in a weighted $L^p$-space, denoted $L^p_{u_\nu}$, with the weight

$$u_\nu(\eta) = |\eta|^\nu e^{|\eta|}, \tag{8}$$

satisfying the following convex conditions:

**1)** $|\psi(\eta)| \leq D|\eta|^z e^{-|\eta|}$, for some $z = z(\eta)$ and $D > 0$;

**2)** $|\psi(\eta) - \psi(\eta - \delta)| \leq \omega_\psi(\eta)|\delta|^\theta$ with $\theta > 0$, $\omega_\psi > 0$, $\omega_\psi \in L^p_{u_\nu}(\mathbb{R})$, $\omega_\psi(\eta) \leq A|\eta|^\gamma e^{-|\eta|}$ with $\gamma = \gamma(\eta)$.

The second set $\mathcal{M}_\mu$ is the set of all $\psi$ in $L^p_{u_\nu}(\mathbb{R})$, such that the following integral is larger than some $\mu > 0$:

$$\mathcal{I}[\psi] = \int_{|\eta| \geq 1} \operatorname{Re}\{\psi(\eta)e^{-i\eta^2}\}|\eta|^\sigma \, d\eta$$

for some $\sigma \in \mathbb{R}$.

The invariance of the compact closure of the convex set $\mathcal{N}_p = \tilde{\mathcal{E}}_D \cap \mathcal{M}_\mu$ under $\mathcal{R}_\beta$ results, via the Tikhonov fixed point theorem, in existence of fixed points $\psi_{\beta,p} \in \overline{\mathcal{N}_p}$ of $\mathcal{R}_\beta$ for every $\beta \in (\beta_0, 1)$. At the same time, as already mentioned, existence of mild blow ups requires existence of one and the same fixed point $\psi_p$ for all $\mathcal{R}_\beta$ with $\beta \in (0,1)$. We show in Section 5 that any limit along a subsequence of the fixed points $\psi_{\beta_i,p}$ with $\beta_i \in (\beta_0, 1)$ is a fixed point of $\mathcal{R}_\beta$ for all $\beta \in (0,1)$. Together with the fact that $\overline{\mathcal{N}}_p \subset L^2_w(\mathbb{R}) \cap L^{\frac{2r+2}{2r+1}}(\mathbb{R})$, where $w(x) = |x|$, this completes the proof of existence of solution whose energy blows up in finite time.

**1.5. Statement of results.** The main result of this paper reads as follows:

**Main theorem.** *For any integer $r \geq 1$ and any $T > 0$ there is a solution $u$ of equation (7) in the class $C^\infty([0,T), C^\infty(\mathbb{R}) \cap L^{2r+2}(\mathbb{R}))$, for which, additionally, $u_x \in C^\infty([0,T), C^\infty(\mathbb{R}) \cap L^2(\mathbb{R}))$, such that the energy (20) become unbounded as $t \to T$.*

We would like to emphasize that although we obtain existence of a blow up solution for the subcritical case $r = 1$, this does not constitute a contradiction with the global possedness results for this case, mentioned in the introduction. Indeed, the initial data for our solutions does not belong to $L^2(\mathbb{R})$ (see Remark 5.3), but rather some higher order norms in $L^{2r+2}(\mathbb{R})$ are bounded for this initial data.

## 2. Prelimenaries

2.1. **The nonlinear Schrödinger equation in Fourier space.** We apply the Fourier transform $v(y) = \mathcal{F}[u](y)$ to (2), so that, using the shorthanded notation $\bullet$ for convolution, the equation (2) can be written as the following integro-differential equation:

$$(9) \quad iv(y,t)_t - y^2 v(y,t) + v \bullet \underbrace{(v \bullet \overline{v \circ (-)}) \bullet \ldots \bullet (v \bullet \overline{v \circ (-)})}_{r \text{ times}} = 0.$$

where $(-)(x) = -x$ and $\overline{v}$ denotes the complex conjugate. We will denote for brevity, the $r$-fold convolution of the functions $v \bullet \overline{v \circ (-)}$ by $(v \bullet \overline{v \circ (-)})^{\bullet r}$.

Multiply (9) by the integrating factor $-ie^{iy^2 t}$ and integrate to get

$$(10) \quad v(y,t) = v(0,y)e^{-iy^2 t} + i \int_0^t \left( v \bullet (v \bullet \overline{v \circ (-)})^{\bullet r} \right)(y,s) e^{-iy^2(t-s)} ds.$$

This is the Picard-Lindelöff equation (7) for the NLS in the Fourier space.

2.2. **Compactness in weighted $L^p$-spaces.** Throughout the paper $L^p_w$ will denote the weighted Lebesgue space with the weight $w$ and the norm

$$(11) \qquad \|f\|_{w,p} := \left( \int_{\mathbb{R}} |f(x)|^p w(x)^p \mathrm{d}x \right)^{\frac{1}{p}}.$$

In what follows we will construct a *convex* $\mathcal{R}_\beta$-invariant, *equitight*, *locally equicontinous* family $\mathcal{N}_p$ of *exponentially decaying* functions in $L^p_{u_\nu}(\mathbb{R})$, $p > 1$. The key result that allows us to claim precompactness of $\mathcal{N}_p$, and, consequently, existence of fixed points, is the following (see [9]):

**Frechet-Kolmogorov-Riesz-Weil Compactness Theorem.** *Let $\mathcal{N}_p$ be a subset of $L^p_w(\mathbb{R})$ with $p \in [1, \infty)$, and let $\tau_\delta f$ denote the translation by $\delta$, $(\tau_\delta f)(x) := f(x - \delta)$. The subset $\mathcal{N}_p$ is relatively compact iff the following properties hold:*

  1) *Equicontinuity:* $\lim_{|\delta| \to 0} \|\tau_\delta f - f\|_{w,p} = 0$ *uniformly in* $\mathcal{N}_p$;
  2) *Equitightness:* $\lim_{r \to \infty} \int_{|x| > r} |f|^p = 0$ *uniformly in* $\mathcal{N}_p$.

By analogy with the renormalization theory in Dynamics, existence of a renormalization-invariant precompact set will be called *apriori* bounds.

## 3. Renormalization problem for a self-similar solution

In this section we will consider an initial value problem for (10) with a self-similar initial condition, and deduce that existence of a self-similar solution to (10) is equivalent to existence of a fixed point for a certain operator.

Consider the equation (10). We will make an assumption that

(12) $$v(y,t) = \tau(t)^{-c}\psi(y\tau(t)),$$

where $\psi$ is some function, and

(13) $$\tau(t) = (T-t)^{\frac{1}{2}}.$$

The following is equation (10) for the initial value problem with the initial condition

(14) $$v_0(y) = \tau_0^{-c}\psi(y\tau_0) = T^{-\frac{c}{2}}\psi(yT^{\frac{1}{2}}), \quad \tau_0 \equiv \tau(0),$$

*under the assumption that $v(y,t)$ remains of the form (12) for $t > 0$*:

$$\frac{\psi(y\tau(t))}{\tau(t)^c} = e^{-iy^2 t}\frac{\psi(y\tau_0)}{\tau_0^c} + i\int_0^t e^{-iy^2(t-s)}\frac{\left(\psi \bullet (\psi \bullet \overline{\psi \circ (-)})^{\bullet r}\right)(y\tau(s))}{\tau(s)^{(2r+1)c+2r}}\,ds,$$

Set $\eta = y\tau_0$, $\xi = y\tau(t)$, $\zeta = y\tau(s)$, then

(15) $$\frac{\psi(\xi)}{\tau(t)^c} = e^{i(\xi^2-\eta^2)}\frac{\psi(\eta)}{\tau_0^c} - 2i\int_\eta^\xi e^{i(\xi^2-\zeta^2)}\frac{\left(\psi \bullet (\psi \bullet \overline{\psi \circ (-)})^{\bullet r}\right)(\zeta)}{y\tau(s)^{(2r+1)c+2r-1}}\,d\zeta.$$

A function $\psi$ solving this equation for all real $\eta$ and $\xi$ such that $|\xi| \leq |\eta|$, provides a solution

(16) $$v(y,t) = (T-t)^{-\frac{c}{2}}\psi\left(y(T-t)^{\frac{1}{2}}\right)$$

to the initial value problem (10) with the initial data (14). By introducing a new time variable,

(17) $$\beta = T^{-\frac{1}{2}}(T-t)^{\frac{1}{2}},$$

equation (15) can be written as

$$\psi(\beta\eta) = \beta^c e^{i(\beta^2\eta^2-\eta^2)}\psi(\eta) - 2i\beta^c|\eta|^c e^{i\beta^2\eta^2}\int_\eta^{\beta\eta} e^{-i\zeta^2}\frac{\left(\psi \bullet (\psi \bullet \overline{\psi \circ (-)})^{\bullet r}\right)(\zeta)}{y|y|^c\tau(s)^{(2r+1)c+2r-1}}\,d\zeta.$$

We can now make the choice

$$1 + c = (2r+1)c + 2r - 1 \implies c = \frac{1-r}{r},$$

then

(18) $$\psi(\beta\eta) = \frac{e^{i(\beta^2\eta^2-\eta^2)}}{\beta^{\frac{r-1}{r}}}\psi(\eta) + \\ + \text{sign}(\eta)\frac{2ie^{i\beta^2\eta^2}}{(\beta|\eta|)^{\frac{r-1}{r}}}\int_{\beta\eta}^\eta e^{-i\zeta^2}\frac{\left(\psi \bullet (\psi \bullet \overline{\psi \circ (-)})^{\bullet r}\right)(\zeta)}{|\zeta|^{\frac{1}{r}}}\,d\zeta,$$

while the problem of existence of solution of (18) can be stated as a fixed point problem for a family of operators

$$\mathcal{R}_\beta[\psi](\eta) = \frac{e^{i\left(1-\frac{1}{\beta^2}\right)\eta^2}}{\beta^{\frac{r-1}{r}}} \psi\left(\frac{\eta}{\beta}\right) +$$

(19)
$$+ \operatorname{sign}(\eta) \frac{2i e^{i\eta^2}}{|\eta|^{\frac{r-1}{r}}} \int_\eta^{\frac{\eta}{\beta}} e^{-i\zeta^2} \frac{\left(\psi \bullet (\psi \bullet \overline{\psi \circ (-)})^{\bullet r}\right)(\zeta)}{|\zeta|^{\frac{1}{r}}} \mathrm{d}\zeta.$$

*Remark* 3.1. It is easy to check that the sets of odd and even functions $\psi$ are invariant under the operator $\mathcal{R}_\beta$.

We emphasize, that to ensure that function $\psi$ generates a self-similar solution (16), we require *one and same* fixed point $\psi$ of $\mathcal{R}_\beta$ for all $\beta \in (0,1]$.

In the next sections we will demonstrate that the operator (19) has a fixed point $\psi_p$ for every $p > 1$, in a certain compact convex subset $\overline{\mathcal{N}}_p \subset L^p_{u_\nu}$ of measurable exponentially decaying functions. The key point is that this $\psi_p$ is fixed by every $\mathcal{R}_\beta$ for all $\beta \in (0,1]$.

In particular, the inverse Fourier transform, $\mathcal{F}^{-1}[\psi_p](x) = \mathcal{F}[\psi](-x)$ is a well defined operator from $\overline{\mathcal{N}}_p$ to, in particular, $L^{2r+2}(\mathbb{R})$, in fact, to $C^\infty(\mathbb{R})$. We will eventually show that the function

$$u(x,t) = (T-t)^{-\frac{1}{2r}} \mathcal{F}^{-1}[\psi_p]\left(\frac{x}{(T-t)^{\frac{1}{2}}}\right)$$

is a $C^\infty(\mathbb{R} \times [0,T))$ solution of (7). For such solution the energy becomes unbounded in finite time. Indeed, assuming for a moment that both $\|\psi_p\|_{w,2}$, $w(x) = |x|$, and $\|\psi_p\|_p$ exist for the specific choice of $p = (2r+2)/(2r+1)$,

$$E(t) = \frac{1}{2}\int_\mathbb{R} |\partial_x u(x,t)|^2 \mathrm{d}x - \frac{1}{2r+2}\int_\mathbb{R} |u(x,t)|^{2r+2} \mathrm{d}x$$

$$= \frac{\tau^{-2(c+1)}}{2}\int_\mathbb{R} \left|\partial_x \mathcal{F}^{-1}[\psi_p]\left(\frac{x}{\tau}\right)\right|^2 \mathrm{d}x - \frac{\tau^{-(c+1)(2r+2)}}{2r+2}\int_\mathbb{R} \left|\mathcal{F}^{-1}[\psi_p]\left(\frac{x}{\tau}\right)\right|^{2r+2} \mathrm{d}x$$

$$= \frac{\tau^{-\frac{2}{r}-1}}{2}\left(\int_\mathbb{R} \left|\partial_\xi \mathcal{F}^{-1}[\psi_p](\xi)\right|^2 \mathrm{d}\xi - \frac{1}{2r+2}\int_\mathbb{R} \left|\mathcal{F}^{-1}[\psi_p](\xi)\right|^{2r+2} \mathrm{d}\xi\right)$$

(20) $$= \tau^{-\frac{2}{r}-1}\left(C_1 \|\psi_p\|_{w,2}^2 + C_2 \|\psi_p\|_p^p\right),$$

where $C_1$ and $C_2$ are some constants and $w(x) = |x|$.

The operator (19) will be referred to as a *renormalization operator*: the equation

$$\psi_p(\eta) = \mathcal{R}_\beta[\psi_p](\eta)$$

says informally that *the time evolution of the initial data at a later time, closer to the blow up time $T$, looked at at a larger scale is equivalent to the initial data itself.*

## 4. A-priori bounds and existence of a renormalization fixed point

### 4.1. A renormalization invariant subset of $L^p_{u_\nu}$-functions.

We will now construct a renormalization invariant subset in $L^p_{u_\nu}(\mathbb{R})$, where $u_\nu$ is as in (8).

**Proposition 4.1.** *For $\forall r \geq 1$, $\forall p > 1$, $\forall \beta_0 \in (0,1)$, there exists $D > 0$, $z_1$ and $z_2$, satisfying*

$$z_1 < -\frac{2r}{2r+1}, \tag{21}$$

$$z_1 > \max\{-\nu - \frac{1}{p}, -1\}, \tag{22}$$

$$z_2 < \min\{-\nu - \frac{1}{p}, -1\}, \tag{23}$$

*such that the subset $\mathcal{E}_D \subset L^p_{u_\nu}(\mathbb{R})$ of exponentially decaying functions,*

$$\mathcal{E}_D = \left\{ f - \text{measurable on } \mathbb{R} : |f(x)| \leq \begin{cases} D|x|^{z_1} e^{-|x|}, & |x| \leq 1, \\ D|x|^{z_2} e^{-|x|}, & |x| > 1 \end{cases} \right\} \tag{24}$$

*is invariant under $\mathcal{R}_\beta$ for all $\beta \in (\beta_0, 1)$.*

*Proof.* First, notice, that for $\mathcal{E}_D$ to be a subset of $L^p_{u_\nu}(\mathbb{R})$, it is sufficient that $z_1 p + \nu p > -1$ and $z_2 p + \nu p < -1$, i.e.

$$z_1 > -\nu - \frac{1}{p}, \tag{25}$$

$$z_z < -\nu - \frac{1}{p}. \tag{26}$$

*1) Linear terms.* As, before, we denote $c = (1-r)/r$.

$$|L(\eta)| \leq \beta^c \left| \psi\left(\frac{\eta}{\beta}\right) \right|$$

$$\leq \begin{cases} \beta^{c-z_1} e^{-(\frac{1}{\beta}-1)|\eta|} D|\eta|^{z_1} e^{-|\eta|}, & |\eta| \leq \beta, \\ \beta^{c-z_2} e^{-(\frac{1}{\beta}-1)|\eta|} D|\eta|^{z_2} e^{-|\eta|}, & |\eta| > \beta, \end{cases}$$

We remark that for $\beta < |\eta| \leq 1$, since $z_1 > z_2$, we have that $|\eta|^{z_2} \leq \beta^{z_2 - z_1} |\eta|^{z_1}$. Therefore,

$$|L(\eta)| \leq \begin{cases} \beta^{c-z_1} e^{-(\frac{1}{\beta}-1)|\eta|} D|\eta|^{z_1} e^{-|\eta|}, & |\eta| \leq 1, \\ \beta^{c-z_2} e^{-(\frac{1}{\beta}-1)|\eta|} D|\eta|^{z_2} e^{-|\eta|}, & |\eta| > 1. \end{cases}$$

$$\leq \begin{cases} (1 + (z_1 - c)\varepsilon) D|\eta|^{z_1} e^{-|\eta|}, & |\eta| \leq 1, \\ (1 + (z_2 - c)\varepsilon) D|\eta|^{z_2} e^{-|\eta|}, & |\eta| > 1. \end{cases} \tag{27}$$

Therefore, the first set of conditions for the invariance of $\mathcal{E}_D$ is $z_1 - c < 0$ and $z_2 - c < 0$.

2) *Nonlinear terms.*

$$|N(\eta)| \leq \frac{2}{\beta^{\frac{r-1}{r}}}|\eta|^{\frac{1-r}{r}} \int_{|\eta|}^{\frac{|\eta|}{\beta}} \frac{(|\psi| \bullet (|\psi| \bullet |\psi \circ (-)|)^{\bullet r})(\zeta)}{|\zeta|^{\frac{1}{r}}} \, d\zeta$$

$$\leq \frac{2}{\beta^{\frac{r-1}{r}}} D^{2r+1}|\eta|^{\frac{1-r}{r}} \int_{|\eta|}^{\frac{|\eta|}{\beta}} \frac{(g \bullet f^{\bullet 2r})(\zeta)}{|\zeta|^{\frac{1}{r}}} d\zeta,$$

where

(28) $$g(x) = \begin{cases} |x|^{k_1} e^{-|x|}, & |x| \leq 1, \\ |x|^{k_2} e^{-|x|}, & |x| > 1 \end{cases}, \quad f(x) = \begin{cases} |x|^{z_1} e^{-|x|}, & |x| \leq 1, \\ |x|^{z_2} e^{-|x|}, & |x| > 1 \end{cases}.$$

We will use a different letter $k$ for the power of the last convolution since we will encounter a case when these powers are indeed distinct later in the paper. Similar integrals will appear several times in this paper, therefore, we will present this computation as a separate lemma. The proof of the lemma can be found in the Appendix 6.

**Lemma 4.2.** *For every $\beta_0 \in (0,1)$, $l > 0$ and every $r \in \mathbb{N}$, $k_1 > -1$, $k_2 < -1$, $z_1 > -1$, $z_2 < -1$, satisfying*

(29) $$z_1 < -\frac{2r-1}{2r},$$

(30) $$2rz_1 + k_1 + 2r < 0,$$

*there exists a constant $K = K(\beta_0, k_1, k_2, z_1, z_2, l, r)$, such that*

$$J_{k_1,k_2,z_1,z_2}^{r,l} := \int_{|\eta|}^{\frac{|\eta|}{\beta}} \frac{(g \bullet f^{\bullet 2r})(\zeta)}{|\zeta|^l} d\zeta$$

(31) $$\leq \varepsilon K |\eta|^n e^{-|\eta|},$$

*where $n = 2rz_1 + k_1 + 2r - l + 1$ for $|\eta| \leq 1$ and $n = \max\{k_2, z_2\} + 1 - l$ for $|\eta| > 1$.*

Since $|N(\eta)| \leq \frac{2}{\beta^{\frac{r-1}{r}}} D^{2r+1}|\eta|^{\frac{1-r}{r}} |J_{z_1,z_2,z_1,z_1}^{r,\frac{1}{r}}(\eta)|$ we get

(32) $$|N(\eta)| \leq \frac{2}{\beta^{\frac{r-1}{r}}} D^{2r+1} \varepsilon K |\eta|^n e^{-|\eta|},$$

where $n = (2r+1)z_1 + 2r$ for $|\eta| \leq 1$, $n = z_2$ for $|\eta| > 1$.

According to (27) and (32), for the set $\mathcal{E}_D$ to be invariant, it is sufficient to require that

(33) $$0 > (z_1 - c)|\eta|^{z_1} + \frac{2}{\beta^{\frac{r-1}{r}}} D^{2r} K |\eta|^{(2r+1)z_1 + 2r},$$

for $|\eta| \leq 1$, and

$$(34) \qquad 0 > (z_2 - c)|\eta|^{z_2} + \frac{2}{\beta^{\frac{r-1}{r}}} D^{2r} |\eta|^{z_2} K.$$

The necessary condition for inequality (33) to hold is

$$(35) \qquad z_1 \leq (2r+1)z_1 + 2r \implies z_1 \geq -1.$$

If (35) holds then (33) can be satisfied by a choice of a sufficiently small $D$.

The second inequality (34) holds if $D$ is chosen sufficienly small. $\square$

The set $\mathcal{E}_D$ contains the trivial function. We will now introduce an extra condition which will define a convex subset of $\mathcal{E}_D$ that does not contain 0.

**Lemma 4.3.** *For every* $D > 0$,

$$(36) \qquad \sigma < -\frac{1}{r}$$

*and* $\beta_0 \in (0,1)$ *there exists* $\mu > 0$, *such that the set*

$$(37) \qquad \mathcal{M}_\mu := \{\psi \in \mathcal{E}_D : \mathcal{I}[\psi] \geq \mu\}, \text{ where}$$

$$\mathcal{I}[\psi] := \int_{|\eta| \geq 1} \mathrm{Re}\,\{\psi(\eta) e^{-i\eta^2}\} |\eta|^\sigma \, d\eta,$$

*is* $\mathcal{R}_\beta$ *invariant for all* $\beta \in (\beta_0, 1)$.

*Proof.* 1) *Linear terms.* Consider $L(\eta)$:

$$\mathcal{I}[L] = \beta^{\frac{1-r}{r}} \int_{|\eta| \geq 1} \mathrm{Re}\left\{e^{-i\frac{\eta^2}{\beta^2}} \psi\left(\frac{\eta}{\beta}\right)\right\} |\eta|^\sigma d\eta$$

$$= \beta^{\sigma + \frac{1}{r}} \int_{|\eta| \geq 1} \mathrm{Re}\left\{e^{-i\frac{\eta^2}{\beta^2}} \psi\left(\frac{\eta}{\beta}\right)\right\} \left|\frac{\eta}{\beta}\right|^\sigma d\frac{\eta}{\beta}$$

$$= \beta^{\sigma + \frac{1}{r}} \int_{|\eta| \geq \frac{1}{\beta}} \mathrm{Re}\left\{e^{-i\eta^2} \psi(\eta)\right\} |\eta|^\sigma \, d\eta$$

$$\geq \beta^{\sigma + \frac{1}{r}} \left(\mathcal{I}[\psi] - 2D\left(\Gamma(z_2 + \sigma + 1, 1) - \Gamma(z_2 + \sigma + 1, \beta^{-1})\right)\right)$$

$$\geq \beta^{\sigma + \frac{1}{r}} \left(\mathcal{I}[\psi] - D\tilde{C}'\varepsilon\right),$$

for some $\tilde{C}'$ depending on $\beta_0$ and $z_2$.

2) *Nonlinear terms.*

$$|\mathcal{I}[N]| = \left|\int_{|\eta|\geq 1} \text{Re}\,\{N(\eta)e^{-i\eta^2}\}|\eta|^\sigma \,\mathrm{d}\eta\right|$$

$$\leq \int_{|\eta|\geq 1} |N(\eta)||\eta|^\sigma \,\mathrm{d}\eta$$

$$\leq \frac{2}{\beta^{\frac{r-1}{r}}} D^{2r+1}\varepsilon \int_{|\eta|\geq 1} K|\eta|^{n+\sigma}e^{-|\eta|}d\eta$$

$$\leq \frac{4}{\beta^{\frac{r-1}{r}}} D^{2r+1}\varepsilon K\Gamma(z_2 + \sigma + 1, 1)$$

$$= D^{2r+1}\varepsilon C''',$$

for some $C'''$ depending on $\beta_0$, $r$, $z_2$ and $z_1$.

We, therefore, obtain that a sufficient condition for the invariance of $\mathcal{R}_\beta$ is

$$\beta^{\sigma+\frac{1}{r}}(\mu - \tilde{C}'\varepsilon D) - D^{2r+1}\varepsilon C'' \geq \mu,$$

which is implied by

(38) $$\left|\sigma + \frac{1}{r}\right|\mu > C'D + C''D^{2r+1},$$

where $C' = \beta_0^{\sigma+\frac{1}{r}}\tilde{C}'$.

The conclusion follows. □

We would now like to show that there exists a choice of constants such that the set $\mathcal{M}_\mu$ is non-empty.

**Lemma 4.4.** *For $\forall r \geq 1$ and $z_1$, $z_2$ satisfying (22), (23), there exist $D > 0$, $\mu > 0$, $\beta_0 \in (0,1)$ and $\sigma$ satisfying (36) such that $\mathcal{M}_\mu$ is non-empty.*

*Proof.* Consider the function

(39) $$\phi(\eta) = \begin{cases} a|\eta|^{z_1}e^{-|\eta|}e^{i\eta^2}, & |\eta| \leq 1, \\ a|\eta|^{z_2}e^{-|\eta|}e^{i\eta^2}, & |\eta| > 1, \end{cases}$$

for some $a < D$. Clearly, $\phi \in \mathcal{E}_D$. We get

$$\mathcal{I}[\phi] = \int_{|\eta|\geq 1} a|\eta|^{z_2+\sigma}e^{-|\eta|}\mathrm{d}\eta.$$

This integral converges, $\mathcal{I}[\phi] = aC''''$, where

$$C'''' = 2\Gamma(z_2 + \sigma + 1, 1).$$

We can now make a choice, for example, $a = D/2$, and

$$\mu = \frac{DC''''}{2}.$$

Then for the condition (38) to be satisfied, it is sufficient that

$$\left|\sigma + \frac{1}{r}\right| C''' \geq 2C' + 2C'' D^{2r},$$

which can be ensured by taking $\sigma$ sufficiently negative, and $D$ sufficiently small. $\square$

4.2. **A-priori bounds.** We start with a straightforward proof of equitightness of the set $\mathcal{E}_M$.

**Proposition 4.5.** *(A-priori bounds: equitightness.) The set $\mathcal{E}_D$ is equitight in $L^p_{u_\mu}$, i. e., for all $\phi \in \mathcal{E}_D$*

$$\lim_{R \to \infty} \int_{|\eta| > R} |\phi(\eta)|^p |\eta|^{\nu p} e^{p|\eta|} d\eta = 0$$

*uniformly on $\mathcal{E}_M$.*

*Proof.*

$$\int_{|\eta|>R} |\phi(\eta)|^p |\eta|^{\nu p} e^{p|\eta|} d\eta \leq D \int_{|\eta|>R} |\eta|^{pz_2 + \nu p} d\eta$$

$$= \frac{2D}{-pz_2 - \nu p - 1} R^{z_2 p + \nu p + 1}.$$

The claim follows since, according to (23), $z_2 p + \nu p + 1 < 0$. $\square$

We can now proceed to equicontinuity.

**Proposition 4.6.** *(A-priori bounds: equicontinuity.) There exists $\delta_0 > 0$, $A > 0$, $D > 0$, $\theta > 0$, $\gamma_2 = z_2 + 1$ and $\gamma_1 = z_1 - \theta$, satisfying*

(40) $$\max\{-\nu - \frac{1}{p}, -1\} < \gamma_1 = z_1 - \theta,$$

(41) $$\gamma_2 = z_2 + 1 < \min\{-\nu - \frac{1}{p}, -1\},$$

*such that the subset*
$\tilde{\mathcal{E}}_D = \{\psi \in \mathcal{E}_D : |\psi(\eta) - \tau_\delta \psi(\eta)| \leq \omega_\psi(\eta)|\delta|^\theta, \omega_\psi(\eta) \leq A|\eta|^{\gamma(|\eta|)} e^{-|\eta|},$

(42) $$\gamma = \gamma_1 \text{ for } |\eta| \leq 1, \gamma = \gamma_2 \text{ for } |\eta| > 1, \omega_\psi \in L^p_{u_\nu}, \ |\delta| < \delta_0\}$$

*of $\mathcal{E}_D$, is renormalization invariant under $\mathcal{R}_\beta$.*

*Proof.* First, we remark that for $\omega_\psi \in L^p_{u_\nu}$, we must have

(43) $$\gamma_1 > -\nu - \frac{1}{p},$$

(44) $$\gamma_2 < -\nu - \frac{1}{p}.$$

We assume that $\psi \in \tilde{\mathcal{E}}_D$ and reproduce these bounds for the renormalized function. We will perform the calculations for $\delta > 0$, the case $\delta < 0$ is similar.

Throughout the proof we will use the notation $C$

$$C = C(z_1, z_2, \gamma_1, \gamma_2, \theta, r, \delta_0, \beta_0)$$

to denote any constant whose specific value is not important.

1) *Linear terms.*

$$\begin{aligned}|L(\eta) - \tau_\delta L(\eta)| &= \beta^c \left| e^{i\left(1-\frac{1}{\beta^2}\right)\eta^2} \psi\left(\frac{\eta}{\beta}\right) - e^{i\left(1-\frac{1}{\beta^2}\right)(\eta-\delta)^2} \psi\left(\frac{\eta-\delta}{\beta}\right) \right| \\ &= \beta^c \left| e^{i\left(1-\frac{1}{\beta^2}\right)\eta^2} - e^{i\left(1-\frac{1}{\beta^2}\right)(\eta-\delta)^2} \right| \left|\psi\left(\frac{\eta}{\beta}\right)\right| + \\ &\quad + \beta^c \left| \psi\left(\frac{\eta}{\beta}\right) - \psi\left(\frac{\eta-\delta}{\beta}\right) \right| \\ &= I_1(\eta) + I_2(\eta).\end{aligned}$$

We consider the first difference:

$$\begin{aligned}I_1(\eta) &= \beta^{c-z} \left| 1 - e^{i\left(1-\frac{1}{\beta^2}\right)(2\eta\delta - \delta^2)} \right| D|\eta|^z e^{-\frac{|\eta|}{\beta}} \\ &= \beta^{c-z} \sqrt{2} \left( 1 - \cos\left( \left(1 - \frac{1}{\beta^2}\right)(2\eta\delta - \delta^2) \right) \right)^{\frac{1}{2}} D|\eta|^z e^{-\frac{|\eta|}{\beta}},\end{aligned}$$

where $z = z_1$ if $|\eta| \leq \beta$ and $z = z_2$ if $|\eta| > \beta$. The difference $1 - \cos x$ can be bounded from above by $x^2/2$ if $|x| \leq 2$ and by $2$ if $|x| > 2$, therefore,

$$(45) \quad I_1(\eta) \leq \beta^{c-z} D \begin{cases} \left(\frac{1}{\beta^2} - 1\right)|2\eta - \delta||\eta|^z e^{-\frac{|\eta|}{\beta}} \delta, & \left|\eta - \frac{\delta}{2}\right| \leq \frac{\beta^2}{\delta(1-\beta^2)}, \\ 2|\eta|^z e^{-\frac{|\eta|}{\beta}}, & \left|\eta - \frac{\delta}{2}\right| > \frac{\beta^2}{\delta(1-\beta^2)}. \end{cases}$$

The difference $I_2$ is straightforward:

$$(46) \quad I_2 \leq \beta^c \omega_\psi\left(\frac{\eta}{\beta}\right) \left|\frac{\delta}{\beta}\right|^\theta \leq A\beta^{c-\gamma-\theta} e^{-\left(\frac{1}{\beta}-1\right)|\eta|} |\eta|^\gamma e^{-|\eta|} \delta^\theta,$$

where $\gamma = \gamma_1$ if $|\eta| \leq \beta$ and $\gamma = \gamma_2$ if $|\eta| > \beta$.

2) *Nonlinear terms. Case* $|\eta| \leq 2\delta$.

$$|N(\eta) - \tau_\delta N(\eta)| \leq 2\left|e^{i\eta^2} - e^{i(\eta-\delta)^2}\right| |\eta|^{\frac{1-r}{r}} \left|\int_\eta^{\frac{\eta}{\beta}} \frac{\psi \bullet (\psi \bullet \overline{\psi \circ (-m)})^{\bullet r}(\zeta)}{|\zeta|^{\frac{1}{r}}} e^{-i\zeta^2} d\zeta\right| +$$

$$+ 2\left||\eta|^{\frac{1-r}{r}} \int_\eta^{\frac{\eta}{\beta}} \frac{\psi \bullet (\psi \bullet \overline{\psi \circ (-m)})^{\bullet r}(\zeta)}{|\zeta|^{\frac{1}{r}}} e^{-i\zeta^2} d\zeta - \right.$$

$$\left. -|\eta-\delta|^{\frac{1-r}{r}} \int_{\eta-\delta}^{\frac{\eta-\delta}{\beta}} \frac{\psi \bullet (\psi \bullet \overline{\psi \circ (-m)})^{\bullet r}(\zeta)}{|\zeta|^{\frac{1}{r}}} e^{-i\zeta^2} d\zeta\right|$$

$$\leq J_1 + J',$$

where

$$J_1 = 2\left|e^{i\eta^2} - e^{i(\eta-\delta)^2}\right| |\eta|^{\frac{1-r}{r}} \left|\int_\eta^{\frac{\eta}{\beta}} \frac{\psi \bullet (\psi \bullet \overline{\psi \circ (-m)})^{\bullet r}(\zeta)}{|\zeta|^{\frac{1}{r}}} e^{-i\zeta^2} d\zeta\right|,$$

$$J' = 2\left||\eta|^{\frac{1-r}{r}} \int_\eta^{\frac{\eta}{\beta}} \frac{\psi \bullet (\psi \bullet \overline{\psi \circ (-m)})^{\bullet r}(\zeta)}{|\zeta|^{\frac{1}{r}}} e^{-i\zeta^2} d\zeta - \right.$$

$$\left. -|\eta-\delta|^{\frac{1-r}{r}} \int_{\eta-\delta}^{\frac{\eta-\delta}{\beta}} \frac{\psi \bullet (\psi \bullet \overline{\psi \circ (-m)})^{\bullet r}(\zeta)}{|\zeta|^{\frac{1}{r}}} e^{-i\zeta^2} d\zeta\right|.$$

We evaluate $J_1$ and $J'$ separately. According to (31),

$$J_1 \leq \sqrt{2}(1 - \cos(2\eta\delta - \delta^2))^{\frac{1}{2}} |\eta|^{\frac{1-r}{r}} 2\varepsilon D^{2r+1} A |\eta|^z e^{-|\eta|}$$

(47) $$\leq \tilde{C}\varepsilon D^{2r+1} |\eta|^z e^{-|\eta|} \begin{cases} |2\eta - \delta|\delta, & \left|\eta - \frac{\delta}{2}\right| \leq \frac{1}{\delta}, \\ 2, & \left|\eta - \frac{\delta}{2}\right| > \frac{1}{\delta}, \end{cases}$$

Since $|\eta| < 2\delta$, we have that

(48) $$J_1 \leq C\varepsilon D^{2r+1} |\eta|^{2r(z_1+1)-1} e^{-|\eta|} \delta^2.$$

Next,

$$J' \leq \tilde{C} D^{2r+1} \varepsilon \left(|\eta|^z e^{-|\eta|} + |\eta-\delta|^z e^{-|\eta-\delta|}\right)$$

(49) $$\leq C D^{2r+1} \varepsilon |\eta|^{z-\theta} e^{-|\eta|} \delta^\theta,$$

where $z = (2r+1)z_1 + 2r$.

3) *Nonlinear terms. Case $|\eta| > 2\delta$*. Throughout this part,

$$z = \begin{cases} (2r+1)z_1 + 2r, & |\eta| \leq 1, \\ z_2, & |\eta| > 1. \end{cases}$$

We start with $J_1$.

$$\begin{aligned} J_1 &\leq \tilde{C}\varepsilon D^{2r+1}|\eta|^z e^{-|\eta|} \begin{cases} |2\eta - \delta|\delta, & \left|\eta - \frac{\delta}{2}\right| \leq \frac{1}{\delta}, \\ 2, & \left|\eta - \frac{\delta}{2}\right| > \frac{1}{\delta}. \end{cases} \\ (50) \qquad &\leq C\varepsilon D^{2r+1} e^{-|\eta|} \begin{cases} \delta|\eta|^{z+1}, & \left|\eta - \frac{\delta}{2}\right| \leq \frac{1}{\delta}, \\ \delta^\theta |\eta|^{z+\theta}, & \left|\eta - \frac{\delta}{2}\right| > \frac{1}{\delta}. \end{cases} \end{aligned}$$

We expand $J'$ as $J' = J_2 + J_3 + J_4 + J_5$ where

$$J_2 = 2|\eta|^{\frac{1-r}{r}} \left| \int_\eta^{\frac{\eta}{\beta}} \frac{\psi \bullet (\psi \bullet \overline{\psi \circ (-m)})^{\bullet r}(\zeta)}{|\zeta|^{\frac{1}{r}}} \left( e^{-i\zeta^2} - e^{-i(\zeta-\delta)^2} \right) d\zeta \right|$$

$$J_3 = 2|\eta|^{\frac{1-r}{r}} \left| \int_\eta^{\frac{\eta}{\beta}} \frac{\psi \bullet (\psi \bullet \overline{\psi \circ (-m)})^{\bullet r}(\zeta) - \psi \bullet (\psi \bullet \overline{\psi \circ (-m)})^{\bullet r}(\zeta - \delta)}{|\zeta|^{\frac{1}{r}}} e^{-i(\zeta-\delta)^2} d\zeta \right|$$

$$J_4 = 2|\eta|^{\frac{1-r}{r}} \left| \int_\eta^{\frac{\eta}{\beta}} \psi \bullet (\psi \bullet \overline{\psi \circ (-m)})^{\bullet r}(\zeta - \delta) \left( \frac{1}{|\zeta - \delta|^{\frac{1}{r}}} - \frac{1}{|\zeta|^{\frac{1}{r}}} \right) e^{-i(\zeta-\delta)^2} d\zeta \right|,$$

$$J_5 = 2 \left| |\eta|^{\frac{1-r}{r}} - |\eta - \delta|^{\frac{1-r}{r}} \right| \left| \int_{\eta-\delta}^{\frac{\eta-\delta}{\beta}} \frac{\psi \bullet (\psi \bullet \overline{\psi \circ (-m)})^{\bullet r}(\zeta)}{|\zeta|^{\frac{1}{r}}} e^{-i\zeta^2} d\zeta \right|,$$

and evaluate $J_i$ separately:

$$J_2 \leq 2|\eta|^{\frac{1-r}{r}} \int_{|\eta|}^{\frac{|\eta|}{\beta}} \frac{|\psi \bullet (\psi \bullet \overline{\psi \circ (-m)})^{\bullet r}(\zeta)|}{|\zeta|^{\frac{1}{r}}} \sqrt{2} \left(1 - \cos(2\zeta\delta - \delta^2)\right)^{\frac{1}{2}} d\zeta$$

$$\leq \begin{cases} 2|\eta|^{\frac{1-r}{r}} \int_{|\eta|}^{\frac{|\eta|}{\beta}} \frac{|\psi\bullet(\psi\bullet\overline{\psi\circ(-m)})^r(\zeta)|}{|\zeta|^{\frac{1}{r}}} |2\zeta - \delta|\delta d\zeta, & \left|\eta - \frac{\delta}{2}\right| \leq \frac{1}{\delta} \\ 4|\eta|^{\frac{1-r}{r}} \int_{|\eta|}^{\frac{|\eta|}{\beta}} \frac{|\psi\bullet(\psi\bullet\overline{\psi\circ(-m)})^r(\zeta)|}{|\zeta|^{\frac{1}{r}}} d\zeta, & \left|\eta - \frac{\delta}{2}\right| > \frac{1}{\delta} \end{cases}$$

$$\leq \begin{cases} 4\delta|\eta|^{\frac{1-r}{r}} \int_{|\eta|}^{\frac{|\eta|}{\beta}} \frac{|\psi\bullet(\psi\bullet\overline{\psi\circ(-m)})^{\frac{1}{r}}(\zeta)|}{|\zeta|^{\frac{1-r}{r}}} d\zeta + \\ \quad + 2\delta^2|\eta|^{\frac{1-r}{r}} \int_{|\eta|}^{\frac{|\eta|}{\beta}} \frac{|\psi\bullet(\psi\bullet\overline{\psi\circ(-m)})^r(\zeta)|}{|\zeta|^{\frac{1}{r}}} d\zeta \\ 4|\eta|^{\frac{1-r}{r}} \int_{|\eta|}^{\frac{|\eta|}{\beta}} \frac{|\psi\bullet(\psi\bullet\overline{\psi\circ(-m)})^r(\zeta)|}{|\zeta|^{\frac{1}{r}}} d\zeta, \end{cases} \begin{array}{l} \beta\left(\frac{\delta}{2} - \frac{1}{\delta}\right) \leq \eta \leq \beta\left(\frac{1}{\delta} + \frac{\delta}{2}\right), \\ \\ \eta < \beta\left(\frac{\delta}{2} - \frac{1}{\delta}\right) \text{ or } \eta > \beta\left(\frac{1}{\delta} + \frac{\delta}{2}\right) \end{array}$$

We conclude, using Lemma 4.2, that

$$J_2 \leq C\delta D^{2r+1}\varepsilon(|\eta|^{z+1} + |\eta|^z \delta)e^{-|\eta|},$$

if

(51) $$\beta\left(\frac{\delta}{2} - \frac{1}{\delta}\right) \leq \eta \leq \beta\left(\frac{1}{\delta} + \frac{\delta}{2}\right).$$

This condition (51) implies that

(52) $$\frac{1}{|\eta|} \geq C\delta$$

for some constant $C$, dependent on $\delta_0$ and $\beta_0$. Therefore, in this case,

(53) $$J_2 \leq CD^{2r+1}\varepsilon(\delta^\theta |\eta|^{z+\theta}e^{-|\eta|} + \delta^2|\eta|^z e^{-|\eta|}).$$

On the other hand, if

$$\eta < \beta\left(\frac{\delta}{2} - \frac{1}{\delta}\right) \text{ or } \eta > \beta\left(\frac{1}{\delta} + \frac{\delta}{2}\right),$$

then $|\eta| \geq C\delta^{-1}$ for some $C$, and

(54) $$\begin{aligned} J_2 &\leq CD^{2r+1}\varepsilon|\eta|^z e^{-|\eta|} \\ &\leq CD^{2r+1}\varepsilon\delta^\theta |\eta|^{z+\theta}e^{-|\eta|}. \end{aligned}$$

Next,

$$J_3 \leq |\eta|^{\frac{1-r}{r}} 2 \int_{|\eta|}^{\frac{|\eta|}{\beta}} \frac{\int_{\mathbb{R}} |(\psi \bullet (\overline{\psi \circ (-)})^{\bullet r})(u)||\psi(u-\zeta) - \psi(u-\zeta+\delta)|\, du}{|\zeta|^{\frac{1}{r}}} d\zeta$$

$$\leq 2|\eta|^{\frac{1-r}{r}} \int_{|\eta|}^{\frac{|\eta|}{\beta}} \frac{\int_{\mathbb{R}} |(\psi \bullet (\overline{\psi \circ (-)})^{\bullet r})(u)| \omega_\psi(u-\zeta) \delta^\theta\, du}{|\zeta|^{\frac{1}{r}}} d\zeta$$

$$\leq CD^{2r} A\delta^\theta |\eta|^{\frac{1-r}{r}} \int_{|\eta|}^{\frac{|\eta|}{\beta}} \frac{(g \bullet f^{\bullet 2r})(\zeta)}{|\zeta|^{\frac{1}{r}}} d\zeta$$

$$\tag{55} \leq CD^{2r} A\varepsilon \delta^\theta |\eta|^\xi e^{-|\eta|},$$

where $g$ and $f$ are as in (28) with $k_1 = \gamma_1$ and $k_2 = \gamma_2$, and

$$\xi = \begin{cases} 2rz_1 + \gamma_1 + 2r, & |\eta| \leq 1, \\ \max\{z_2, \gamma_2\}, & |\eta| > 1. \end{cases}$$

To compute $J_4$, we notice that

$$|\zeta| > 2\delta \implies \frac{|\zeta - \delta|}{|\zeta|} < \frac{3}{2} \implies$$

$$\left| \frac{1}{|\zeta|^{\frac{1}{r}}} - \frac{1}{|\zeta - \delta|^{\frac{1}{r}}} \right| = \frac{\left| |\zeta - \delta|^{\frac{1}{r}} - |\zeta|^{\frac{1}{r}} \right|}{|\zeta|^{\frac{1}{r}} |\zeta - \delta|^{\frac{1}{r}}}$$

$$\leq \frac{1}{r} \frac{\delta}{|\zeta||\zeta - \delta|^{\frac{1}{r}}}$$

$$\leq \frac{3}{2r} \frac{\delta}{|\zeta - \delta|^{\frac{1}{r}+1}}$$

$$\leq \frac{3}{2r} \frac{\delta^\theta}{|\zeta - \delta|^{\frac{1}{r}+\theta}}.$$

Therefore,

$$\tag{56} J_4 \leq C\varepsilon D^{2r+1} \delta^\theta |\eta|^{z+\theta} e^{-|\eta|}.$$

Finally, to compute $J_5$, we consider $||\eta|^c - |\eta - \delta|^c|$:

$$\left| |\eta|^{\frac{1-r}{r}} - |\eta - \delta|^{\frac{1-r}{r}} \right| \leq |\eta|^{\frac{1-r}{r}} \left( \left(1 - \frac{\delta}{|\eta|}\right)^{\frac{1-r}{r}} - 1 \right)$$

$$\leq |\eta|^{\frac{1-r}{r}} 2^{\frac{r-1}{r}+1} \frac{r-1}{r} \frac{\delta}{|\eta|}.$$

Therefore,

$$J_5 \leq CD^{2r+1}\delta\varepsilon|\eta|^{z-1}e^{-|\eta|}$$

(57)
$$\leq CD^{2r+1}\begin{cases} \delta^\theta\varepsilon|\eta|^{z-\theta}e^{-|\eta|}, & 2\delta \leq |\eta| \leq 1, \\ \delta\varepsilon|\eta|^{z-1}e^{-|\eta|}, & |\eta| > 1. \end{cases}$$

*4) Collecting all terms.* We start by remarking that if $c - \gamma - \theta > 0$, then, according to (46),

$$I_2 \leq A(1 + C(\gamma + \theta - c)\varepsilon)|\eta|^\gamma e^{-|\eta|}\delta^\theta$$

for all $\eta$, where $\gamma = \gamma_1$ if $|\eta| \leq \beta$ and $\gamma = \gamma_2$ if $|\eta| > \beta$. Furthermore, according to (45),

$$I_1(\eta) \leq DCe^{-|\eta|}\begin{cases} \varepsilon\delta^2|\eta|^z, & |\eta| \leq \delta \text{ and } \left|\eta - \frac{\delta}{2}\right| \leq \frac{\beta^2}{\delta(1-\beta^2)} \\ \varepsilon\delta|\eta|^{z+1}, & |\eta| > \delta \text{ and } \left|\eta - \frac{\delta}{2}\right| \leq \frac{\beta^2}{\delta(1-\beta^2)} \\ |\eta|^z, & \left|\eta - \frac{\delta}{2}\right| > \frac{\beta^2}{\delta(1-\beta^2)}. \end{cases}$$

$$\leq DCe^{-|\eta|}\begin{cases} \varepsilon\delta^2|\eta|^z, & |\eta| \leq \delta \text{ and } \left|\eta - \frac{\delta}{2}\right| \leq \frac{\beta^2}{\delta(1-\beta^2)} \\ \varepsilon\delta|\eta|^{z+1}, & |\eta| > \delta \text{ and } \left|\eta - \frac{\delta}{2}\right| \leq \frac{\beta^2}{\delta(1-\beta^2)} \\ \varepsilon\delta|\eta|^{z+1}, & \left|\eta - \frac{\delta}{2}\right| > \frac{\beta^2}{\delta(1-\beta^2)}. \end{cases}$$

$$\leq DC\varepsilon\delta e^{-|\eta|}\begin{cases} |\eta|^z, & |\eta| \leq 1 \\ |\eta|^{z+1}, & |\eta| > 1 \end{cases}.$$

According to (47), (50), (49), (53), (54), (55), (56) and (57), the nonlinear terms can be bounded as

(58) $$|N(\eta) - \tau_\delta N(\eta)| \leq C\varepsilon D^{2r}\delta^\theta(D|z|^\zeta + A|z|^\xi)e^{-|\eta|},$$

where $\zeta = z_1 - \theta$ if $|\eta| \leq 1$, $\zeta = z_2 + \theta$ if $|\eta| > 1$, and $\xi = 2rz_1 + \gamma_1 + 2r$ if $|\eta| \leq 1$, $\xi = \max\{z_2, \gamma_2\}$ if $|\eta| > 1$. Therefore, for the invariance it is sufficient ot choose $\gamma = z_1 - \theta$ if $|\eta| \leq 1$ and $\gamma = z_2 + 1$ if $|\eta| > 1$.

The invariance equation now becomes

$$A(1 + C(\gamma + \theta - c)\varepsilon) + CD\varepsilon + CD^{2r}\varepsilon(D + A) \leq A,$$

which can be insured, if $\gamma + \theta - c < 0$ by a choice of a sufficiently large $A$ and sufficiently small $D$. □

It is easy to check that the function (39) is in $\tilde{\mathcal{E}}_D$, therfore

**Lemma 4.7.** *The intersection of the set $\mathcal{M}_\mu$ from Lemma 4.4 with $\tilde{\mathcal{E}}_D$ is non-empty.*

The Porpositions 4.5 and 4.6 about the renormalization invariance of an equitight and equicontinous set in $L^p_{u_\nu}$ will be refered to as *apriori bounds*. Apriori bounds result in existence of a renormalization fixed point through the following straightworfard argument.

**Corollary 4.8.** *For every $p > 1$, there exists $\nu$ and $\beta_0 \in (0,1)$, such that the operator $\mathcal{R}_\beta$ has a non-trivial fixed point $\psi_{\beta,p}$ in $L^p_{u_\nu}$ for all $\beta \in (\beta_0, 1)$.*

*Proof.* The intersection $\mathcal{N}_p$ of the set $\tilde{\mathcal{E}}_D$, described in Proposition 4.6, and the set $\mathcal{M}_\mu$, described in Lemma 4.3, is non-empty, accroding to Lemma 4.7, convex, as the intersection of two convex sets, pre-compact, by the Frechet-Kolmogorov-Riesz-Weil Compactness Theorem 2.2, and, for every $\beta \in (\beta_0, 1)$ (for some $\beta_0 \in (0,1)$), is $\mathcal{R}_\beta$-invariant, by Lemma 4.3 and Propositions 4.5 and 4.6.

The Tychonoff Fixed Point Theorem implies existence of a fixed point $\psi_{\beta,p} \in \overline{\mathcal{N}_p} \subset L^p_{u_\nu}$ for the continous operator $\mathcal{R}_\beta$. □

## 5. Limits of renormalization fixed points

We emphasize, that existence of fixed points for the family $\mathcal{R}_\beta$ does not yet imply existence of a solution to the equation (15). Rather, one needs to show, that there is *one and the same* fixed point $\psi_p$ for *all* $\beta \in (0,1)$. The next Lemma provides a step in this direction.

**Lemma 5.1.** *Any fixed point $\psi_\beta$ of the operator $\mathcal{R}_\beta$ is also fixed by $\mathcal{R}_{\beta^n}$ for all $n \in \mathbb{N}$.*

*Proof.* We will use the following form of the fixed point equation $\mathcal{R}_\beta[\psi_\beta] = \psi_\beta$:

$$\beta^{\frac{r-1}{r}} e^{-i\beta^2 \eta^2} \psi_\beta(\beta \eta) = e^{-i\eta^2} \psi_\beta(\eta) +$$

$$(59) \qquad + 2i \operatorname{sign}(\eta) |\eta|^{\frac{1-r}{r}} \int_{\beta\eta}^{\eta} e^{-i\zeta^2} \frac{\psi_\beta \bullet (\psi_\beta \bullet \overline{\psi_\beta \circ (-)})^{\bullet r}(\zeta)}{|\zeta|^{\frac{1}{r}}} d\zeta,$$

For notational brevity we will denote the integral appearing in thes equation as $K_\beta(\eta)$:

$$K_\beta(\eta) = 2i \operatorname{sign}(\eta) |\eta|^{\frac{1-r}{r}} \int_{\beta\eta}^{\eta} e^{-i\zeta^2} \frac{\psi_\beta \bullet (\psi_\beta \bullet \overline{\psi_\beta \circ (-)})^{\bullet r}(\zeta)}{|\zeta|^{\frac{1}{r}}} d\zeta,$$

and we will repeatedly use the fact that

$$K_{\beta^k}(x) + \beta^{-ck} K_\beta(\beta^k x) = K_{\beta^{k+1}}(\eta).$$

We have at the base of induction, evaluating (59) at $\eta = \beta x$:

$$\beta^{-c} e^{-i\beta^4 x^2} \psi_\beta(\beta^2 x) = e^{-i\beta^2 x^2} \psi_\beta(\beta x) + K_\beta(\beta x)$$
$$= \beta^c e^{-i\eta^2} \psi_\beta(x) + \beta^c K_\beta(x) + K_\beta(\beta x)$$
$$= \beta^c \left( e^{-i\eta^2} \psi_\beta(x) + K_{\beta^2}(x) \right),$$

i.e.

$$\beta^{-2c} e^{-i\beta^4 x^2} \psi_\beta(\beta^2 x) = e^{-i\eta^2} \psi_\beta(x) + K_{\beta^2}(x),$$

and $\psi_\beta$ is a fixed point for $\mathcal{R}_{\beta^2}$.

Assume the result for all $n \leq k$ and evaluate (59) at $\eta = \beta^k x$:

$$\beta^{-c} e^{-i\beta^{2(k+1)}x^2} \psi_\beta(\beta^{k+1}x) = e^{-i\beta^{2k}x^2}\phi_\beta(\beta^k x) + K_\beta(\beta^k x)$$
$$= \beta^{ck} e^{-ix^2}\phi_\beta(x) + \beta^{ck} K_{\beta^k}(x) + K_\beta(\beta^k x)$$
$$= \beta^{ck}\left(e^{-ix^2}\phi_\beta(x) + K_{\beta^{k+1}}(x)\right),$$

i.e.

$$\beta^{-c(k+1)} e^{-i\beta^{2(k+1)}x^2}\psi_\beta(\beta^{k+1}x) = e^{-ix^2}\phi_\beta(x) + K_{\beta^{k+1}}(x).$$

$\square$

We would like now to address the question of what happens to the functions $\psi_{\beta,p}$, described in Corollary 4.8, as $\beta$ approaches 1.

Consider a sequence of $\{\psi_{\beta_i,p}\}_{i=0}^\infty$, $\beta_i \in (\beta_0, 1)$ of fixed points of operators $\mathcal{R}_{\beta_i}$ as $\beta_i \to 1$ with $\psi_{\beta_i,p} \in \overline{\mathcal{N}_p}$ where $\overline{\mathcal{N}_p}$ is as in Proposition 4.8. By compactness of the set $\overline{\mathcal{N}_p}$, we can pass to a converging subsequence, which we will also denote $\{\psi_{\beta_i,p}\}_{i=0}^\infty$. Its limit $\psi_p$ is non-trivial.

Consider this sequence $\{\psi_{\beta_i,p}\}_{i=0}^\infty$, $L^p_{u_\nu}$-convering to $\psi_p$. For any $\beta \in (0,1)$ there exists a diverging sequence of integers $n_i$, such that $\beta_i^{n_i} \to \beta$. Indeed, for any $\varepsilon > 0$, a sufficient condition for $\beta_i^{n_i} \in (\beta - \varepsilon, \beta + \varepsilon)$ is

$$n_i \in \left(\frac{\ln(\beta+\varepsilon)}{\ln \beta_i}, \frac{\ln(\beta-\varepsilon)}{\ln \beta_i}\right).$$

The length of this interval is of the order $-\varepsilon/\beta \ln \beta_i$. Therefore, for any $\varepsilon$ there exists $I \in \mathbb{N}$, such that this length is larger than 1 for all $i \geq I$, and the interval contains some positive integer $n_i$. Therefore,

$$\|\mathcal{R}_\beta[\psi_{\beta_i,p}] - \psi_{\beta_i,p}\|_{u_\nu,p} \leq \|(\mathcal{R}_{\beta_i^{n_i}}[\psi_{\beta_i,p}] - \psi_{\beta_i,p}) + (\mathcal{R}_\beta[\psi_{\beta_i,p}] - \mathcal{R}_{\beta_i^{n_i}}[\psi_{\beta_i,p}])\|_{u_\nu,p}$$
$$\leq \|\mathcal{R}_{\beta_i^{n_i}}[\psi_{\beta_i,p}] - \psi_{\beta_i,p}\|_{u_\nu,p} +$$
$$+ \beta^c \left\|e^{i\left(1-\frac{1}{\beta^2}\right)\eta^2}\psi_{\beta_i,p} \circ \beta^{-1} - e^{i\left(1-\frac{1}{\beta_i^{2n_i}}\right)\eta^2}\psi_{\beta_i,p} \circ \beta_i^{-n_i}\right\|_{u_\nu,p} +$$
$$+ |\beta^c - \beta_i^{cn_i}| \left\|e^{i\left(1-\frac{1}{\beta_i^{2n_i}}\right)\eta^2}\psi_{\beta_i,p} \circ \beta_i^{-n_i}\right\|_{u_\nu,p} +$$
$$+ \left\|2e^{i\eta^2}|\eta|^c \int_{\frac{\eta}{\beta}}^{\frac{\eta}{\beta_i^{n_i}}} e^{-i\zeta^2}\frac{\psi_{\beta_i,p} \bullet (\psi_{\beta_i,p} \bullet \overline{\psi_{\beta_i,p} \circ (-)})^{\bullet r}(\zeta)}{|\zeta|^{\frac{1}{r}}} d\zeta\right\|_{u_\nu,p}.$$

According to Lemma (5.1), we have $\mathcal{R}_{\beta_i^{n_i}}[\psi_{\beta_i,p}] = \psi_{\beta_i,p}$, therefore,

$$\|\mathcal{R}_\beta[\psi_{\beta_i,p}] - \psi_{\beta_i,p}\|_{u_\nu,p} \leq \beta^c \left\| e^{i\left(1-\frac{1}{\beta^2}\right)\eta^2} \psi_{\beta_i,p} \circ \beta^{-1} - e^{i\left(1-\frac{1}{\beta_i^{2n_i}}\right)\eta^2} \psi_{\beta_i,p} \circ \beta^{-1} \right\|_{u_\nu,p} +$$

$$+ |\beta^c - \beta_i^{cn_i}| \left\| e^{i\left(1-\frac{1}{\beta_i^{2n_i}}\right)\eta^2} \psi_{\beta_i,p} \circ \beta_i^{-n_i} \right\|_{u_\nu,p} +$$

$$+ \beta^c \left\| e^{i\left(1-\frac{1}{\beta_i^{2n_i}}\right)\eta^2} \psi_{\beta_i,p} \circ \beta^{-1} - e^{i\left(1-\frac{1}{\beta_i^{2n_i}}\right)\eta^2} \psi_{\beta_i,p} \circ \beta_i^{-n_i} \right\|_{u_\nu,p} +$$

$$+ \left\| 2e^{i\eta^2} |\eta|^c \int_{\frac{\eta}{\beta}}^{\frac{\eta}{\beta_i^{n_i}}} e^{-i\zeta^2} \frac{\psi_{\beta_i,p} \bullet (\psi_{\beta_i,p} \bullet \overline{\psi_{\beta_i,p} \circ (-)})^{\bullet r}(\zeta)}{|\zeta|^{\frac{1}{r}}} \, d\zeta \right\|_{u_\nu,p}$$

$$(60) \qquad \leq C\sqrt{2} \left\| \left(1 - \cos\left(\left(\frac{1}{\beta_i^{2n_i}} - \frac{1}{\beta^2}\right)\eta^2\right)\right) \psi_{\beta_i,p} \circ \beta^{-1} \right\|_{u_\nu,p} +$$

$$+ C|\beta - \beta_i^{n_i}| \left\| \psi_{\beta_i,p} \circ \beta_i^{-n_i} \right\|_{u_\nu,p} +$$

$$(61) \qquad + \beta^c \left\| e^{\left(1-\frac{1}{\beta_i^{2n_i}}\right)\eta^2} \psi_{\beta_i,p} \circ \beta^{-1} - e^{\left(1-\frac{1}{\beta_i^{2n_i}}\right)\eta^2} \psi_{\beta_i,p} \circ \beta_i^{-n_i} \right\|_{u_\nu,p} +$$

$$(62) \qquad + \left\| 2e^{i\eta^2} |\eta|^c \int_{\frac{\eta}{\beta}}^{\frac{\eta}{\beta_i^{n_i}}} e^{-i\zeta^2} \frac{\psi_{\beta_i,p} \bullet (\psi_{\beta_i,p} \bullet \overline{\psi_{\beta_i,p} \circ (-)})^{\bullet r}(\zeta)}{|\zeta|^{\frac{1}{r}}} \, d\zeta \right\|_{u_\nu,p}.$$

where $C$ is some constant.

The first norm (60) has been already considered in Proposition 4.6, and is bounded by $C|\beta - \beta_i^{n_i}|$.

The last norm (62) can be bound using (32) with $\varepsilon$ substituted by $\beta - \beta_i^{n_i}$ and $n, m, k$ equal to either $z_1$ or $z_2$. In particular the condition of convergence of this integral is

$$(63) \qquad z_1 > -\nu - \frac{1}{p}.$$

The remaining norm (61) can be split in three parts: $I_1 + I_2 + I_3$, where

$$I_1^p = \beta^{cp} \int_{|\eta| < \frac{\delta_0}{\left|\frac{1}{\beta} - \frac{1}{\beta_i^{n_i}}\right|}} \left|\psi_{\beta_i,p}\left(\frac{\eta}{\beta}\right) - \psi_{\beta_i,p}\left(\frac{\eta}{\beta_i^{n_i}}\right)\right|^p |\eta|^{\nu p} e^{p|\eta|} d\eta,$$

$$I_2^p = \beta^{cp} \int_{|\eta| \geq \frac{\delta_0}{\left|\frac{1}{\beta} - \frac{1}{\beta_i^{n_i}}\right|}} \left|\psi_{\beta_i,p}\left(\frac{\eta}{\beta}\right)\right|^p |\eta|^{\nu p} e^{p|\eta|} d\eta,$$

$$I_3^p = \beta^{cp} \int_{|\eta| \geq \frac{\delta_0}{\left|\frac{1}{\beta} - \frac{1}{\beta_i^{n_i}}\right|}} \left|\psi_{\beta_i,p}\left(\frac{\eta}{\beta_i^{n_i}}\right)\right|^p |\eta|^{\nu p} e^{p|\eta|} d\eta.$$

The difference of the arguments of $\psi_{\beta_i,p}$ in $I_1$ is less than $\delta_0$, and we can use Proposition (4.6), therefore these integrals can be estimated as follows:

$$I_1^p \leq \beta^{cp} \int_{|\eta| < \frac{\delta_0}{\left|\frac{1}{\beta} - \frac{1}{\beta_i^{n_i}}\right|}} \omega_{\psi_{\beta_i,p}}^p\left(\frac{\eta}{\beta}\right) \left|\frac{\eta}{\beta} - \frac{\eta}{\beta_i^{n_i}}\right|^{p\theta} |\eta|^{\nu p} e^{p|\eta|} d\eta$$

$$\leq C|\beta - \beta_i^{n_i}|^{p\theta},$$

$$I_k^p \leq C \int_{|\eta| \geq \frac{\delta_0}{\left|\frac{1}{\beta} - \frac{1}{\beta_i^{n_i}}\right|}} |\eta|^{z_2 p + \nu p} d\eta$$

$$\leq C |\beta - \beta_i^{n_i}|^{-1-z_2 p - \nu p},$$

$k = 2, 3$, where we have bounded the remaining integrals with constants, using the fact that $\omega_\psi$ and $\psi_{\beta_i,p}$ are exponentially bounded (here, since their specific values are irrelevant, all constants are denoted with $C$). Notice that the power $-1 - z_2 p - \nu p$ is positive by (23).

We obtain that

(64) $$\lim_{i \to \infty} \|\mathcal{R}_\beta[\psi_{\beta_i,p}] - \psi_{\beta_i,p}\|_{u_\nu,p} = 0.$$

The immediate consequence of (64) is

**Proposition 5.2.** *For every $r \in \mathbb{N}$, $r \geq 1$, $p > 1$ and $\nu \in \mathbb{R}$, satisfying*

$$\nu + \frac{1}{p} > \frac{2r}{2r+1},$$

*there exists $\psi_p \in \overline{\mathcal{N}}_p \subset L_{u_\nu}^p(\mathbb{R})$, which is a fixed point of $\mathcal{R}_\beta$ for all $\beta \in (0,1)$.*

Recall that for any $p \in (1, 2]$ the Fourier transform is a bounded operator from $L^p$ to $L^{p'}$, where $p'$ is the Hölder conjugate of $p$. Set $p' = 2r + 2$ and

$\nu = 0$, then the conditions (21), (22), (23) become

(65) $$z_2 < -1 < -\frac{2r+1}{2r+2} < z_1 < -\frac{2r}{2r+1},$$

therefore, for any integre $r \geq 1$, there is a choice of $z_1$ and $z_2$ that satisfy the inequalities (65).

Such choice of $p = (2r+2)/(2r+1)$, $\nu = 0$, $z_1$ and $z_2$ defines a relatively compact set $\mathcal{N}_p^*$.

Furthermore, notice that if we set $\nu = 1$ and $p = 2$, then we have that for any $r \geq 0$,

$$-\nu - \frac{1}{p} = -\frac{3}{2} < z_1 < -\frac{2r}{2r+1}$$

which means that the $\mathcal{N}_p^* \subset L_w^2(\mathbb{R})$ with $w(x) = |x|$.

We conclude that $\psi_p \in L_w^2(\mathbb{R}) \cap L^{\frac{2r+2}{2r+1}}(\mathbb{R})$, while

(66) $$\mathcal{F}^{-1}[\psi_p] \in C^\infty(\mathbb{R}) \cap L^{2r+2}(\mathbb{R}), \quad \partial \mathcal{F}^{-1}[\psi_p] \in L^2(\mathbb{R}).$$

*Remark* 5.3. A crucial remark is that we can not claim that $\mathcal{F}^{-1}[\psi_p] \in H^1(\mathbb{R})$, the standard Sobolev space. Indeed, choosing $\nu = 0$ and $p = 2$ in the condition $z_1 > -\nu - 1/p$, we get that $r$ must satisfy

$$-\frac{1}{2} < z_1 < -\frac{2r}{2r+1} \implies r < \frac{1}{2}.$$

Since $\mathcal{F}^{-1}[\psi_p]$ is in $L^{2r+2}(\mathbb{R})$ and $\partial \mathcal{F}^{-1}[\psi_p] \in L^2(\mathbb{R})$, the energy (20) is well defined and bounded for all $t < T$, becoming unbounded as $t$ approaches $T$:

$$E(t) = \frac{1}{2}\int_\mathbb{R} |\partial_x u(x,t)|^2 dx - \frac{1}{2r+2}\int_\mathbb{R} |u(x,t)|^{2r+2} dx$$
$$= \frac{\tau^{-\frac{2}{r}-1}}{2}\left(\int_\mathbb{R} |\partial_\xi \mathcal{F}^{-1}[\psi_p](\xi)|^2 d\xi - \frac{1}{r+1}\int_\mathbb{R} |\mathcal{F}^{-1}[\psi_p](\xi)|^{2r+2} d\xi\right)$$
$$= \tau^{-\frac{2}{r}-1}\left(C_1\|\psi_p\|_w^2 + C_2\|\psi_p\|_p^p\right),$$

where $C_1$ and $C_2$ are some constants and $w(x) = |x|$.

We conclude with the following

**Theorem 5.4.** *For any integer $r \geq 1$ the equation (7) has a solution*

$$u(x,t) = (T-t)^{\frac{r-1}{2r}}\mathcal{F}^{-1}\left[\psi_p\left((\cdot)(T-t)^{\frac{1}{2}}\right)\right](x)$$
$$= (T-t)^{-\frac{1}{2r}}\mathcal{F}^{-1}[\psi_p]\left(\frac{x}{(T-t)^{\frac{1}{2}}}\right),$$

*which belongs to the class $C^\infty\left([0,T), C^\infty(\mathbb{R}) \cap L^{2r+2}(\mathbb{R})\right)$, while $u_x \in C^\infty([0,T), C^\infty(\mathbb{R}) \cap L^2(\mathbb{R}))$.*

*Additionally, the energy (20) becomes unbounded in finite time.*

## 6. Appendix: Proof of Lemma 4.2

We first consider one convolution $g \bullet f$ where $g$ and $f$ are as in (28),

$$(g \bullet f)(u) := \int_{-\infty}^{\infty} |y|^k |u-y|^z e^{-|y|-|u-y|} dy$$

$$= \int_{-\infty}^{\infty} |u|^{k+z+1} |x|^k |1-x|^z e^{-|u||x|-|u||1-x|} dy$$

$$= \int_{x \leq 0} |u|^{k+z+1} |x|^k (1-x)^z e^{2|u|x-|u|} dx + \int_{0 < x \leq 1} |u|^{k+z+1} x^k (1-x)^z e^{-|u|} dx +$$

$$+ \int_{x > 1} |u|^{k+z+1} x^k (x-1)^z e^{|u|-2|u|x} dx$$

$$=: I_1 + I_2 + I_3,$$

where we emphasize that the powers $k$ and $z$ are functions of $|u||x|$ and $|u||1-x|$ respectively. We have

$$I_1 = \int_{\substack{x \leq 0, |u||x| \leq 1 \\ |u|(1-x) \leq 1}} |u|^{k_1+z_1+1} |x|^{k_1} (1-x)^{z_1} e^{2|u|x-|u|} dx + \int_{\substack{x \leq 0, |u||x| \leq 1 \\ |u|(1-x) > 1}} |u|^{k_1+z_2+1} |x|^{k_1} (1-x)^{z_2} e^{2|u|x-|u|} dx +$$

$$+ \int_{\substack{x \leq 0, |u||x| > 1 \\ |u|(1-x) \leq 1}} |u|^{k_2+z_1+1} |x|^{k_2} (1-x)^{z_1} e^{2|u|x-|u|} dx + \int_{\substack{x \leq 0, |u||x| > 1 \\ |u|(1-x) > 1}} |u|^{k_2+z_2+1} |x|^{k_2} (1-x)^{z_2} e^{2|u|x-|u|} dx$$

$$= I_1^1 + I_1^2 + I_1^3 + I_1^4,$$

and

$$I_2 = \int_{\substack{0 < x \leq 1, |u||x| \leq 1 \\ |u|(1-x) \leq 1}} |u|^{k_1+z_1+1} |x|^{k_1} (|x|-1)^{z_1} e^{-|u|} dx + \int_{\substack{0 < x \leq 1, |u||x| \leq 1 \\ |u|(1-x) > 1}} |u|^{k_1+z_2+1} |x|^{k_1} (1-x)^{z_2} e^{-|u|} dx +$$

$$+ \int_{\substack{0 < x \leq 1, |u||x| > 1 \\ |u|(1-x) \leq 1}} |u|^{k_2+z_1+1} |x|^{k_2} (1-x)^{z_1} e^{-|u|} dx + \int_{\substack{0 < x \leq 1, |u||x| > 1 \\ |u|(1-x) > 1}} |u|^{k_2+z_2+1} |x|^{k_2} (1-x)^{z_2} e^{-|u|} dx$$

$$= I_2^1 + I_2^2 + I_2^3 + I_2^4,$$

and
$$I_3 = \int_{\substack{x>1,|u||x|\leq 1 \\ |u|(x-1)\leq 1}} |u|^{k_1+z_1+1}|x|^{k_1}(x-1)^{z_1}e^{|u|-2|u|x}\mathrm{d}x + \int_{\substack{x>1,|u||x|\leq 1 \\ |u|(x-1)>1}} |u|^{k_1+z_2+1}|x|^{k_1}(x-1)^{z_2}e^{|u|-2|u|x}\mathrm{d}x +$$
$$+ \int_{\substack{x>1,|u||x|>1 \\ |u|(x-1)\leq 1}} |u|^{k_2+z_1+1}|x|^{k_2}(x-1)^{z_1}e^{|u|-2|u|x}\mathrm{d}x + \int_{\substack{x>1,|u||x|>1 \\ |u|(x-1)>1}} |u|^{k_2+z_2+1}|x|^{k_2}(x-1)^{z_2}e^{|u|-2|u|x}\mathrm{d}x$$
$$= I_3^1 + I_3^2 + I_3^3 + I_3^4.$$

We evaluate these integrals in three cases separately, writing out only the integrals with non-empty ranges of integration.

1) Case $0 \leq |u| \leq 1$.

$$I_1 = I_1^1 + I_1^2 + I_1^4$$
$$= \int_{1-\frac{1}{|u|}\leq x\leq 0} |u|^{k_1+z_1+1}|x|^{k_1}(1-x)^{z_1}e^{2|u|x-|u|}\mathrm{d}x + \int_{-\frac{1}{|u|}\leq x<1-\frac{1}{|u|}} |u|^{k_1+z_2+1}|x|^{k_1}(1-x)^{z_2}e^{2|u|x-|u|}\mathrm{d}x +$$
$$+ \int_{x<-\frac{1}{|u|}} |u|^{k_2+z_2+1}|x|^{k_2}(1-x)^{z_2}e^{2|u|x-|u|}\mathrm{d}x$$
$$\leq \int_1^{\frac{1}{|u|}}(x-1)^{k_1}x^{z_1}\mathrm{d}x e^{-|u|}|u|^{z_1+k_1+1} + \int_{2(1-|u|)}^2 x^{k_1}e^{-x}\mathrm{d}x\frac{1}{2^{k_1+1}}e^{-|u|} +$$
$$+ \int_2^\infty x^{k_2}e^{-x}\mathrm{d}x\frac{(1+|u|)^{z_2}}{2^{k_2+1}}e^{-|u|}$$
$$\leq \int_1^\infty (x-1)^{k_1}x^{z_1}\mathrm{d}x e^{-|u|}|u|^{z_1+k_1+1} + \int_{2(1-|u|)}^2 x^{k_1}e^{-x}\mathrm{d}x\frac{1}{2^{k_1+1}}e^{-|u|} +$$
$$+ \int_2^\infty x^{k_2}e^{-x}\mathrm{d}x\frac{(1+|u|)^{z_2}}{2^{k_2+1}}e^{-|u|}.$$

The condition for convergence of the first integral is

(67) $$z_1 + k_1 + 1 < 0,$$
(68) $$k_1 > -1.$$

Then
$$I_1 \leq B(-1-k_1-z_1, k_1+1)|u|^{z_1+k_1+1}e^{-|u|} +$$
$$+ (\gamma(k_1+1, 2) - \gamma(k_1+1, 2(1-|u|)))\frac{e^{-|u|}}{2^{k_1+1}} +$$
$$+ \Gamma(k_2+1, 2)\frac{(1+|u|)^{z_2}}{2^{k_2+1}}e^{-|u|},$$

The next set of integrals is as follows:

$$I_2 = I_2^1$$
$$= \int_{0<x\leq 1} |u|^{k_1+z_1+1}|x|^{k_1}(1-x)^{z_1}e^{-|u|}dx.$$

The condition for convergence of this integral is

$$k_1 > -1, \quad z_1 > -1.$$

Then,

$$I_2 \leq B(k_1+1, z_1+1)|u|^{k_1+z_1+1}e^{-|u|}.$$

Finally,

$$I_3 = I_3^1 + I_3^3 + I_3^4$$
$$= \int_{1<x\leq \frac{1}{|u|}} |u|^{k_1+z_1+1}|x|^{k_1}(x-1)^{z_1}e^{|u|-2|u|x}dx + \int_{\frac{1}{|u|}<x\leq \frac{1}{|u|}+1} |u|^{k_2+z_1+1}|x|^{k_2}(x-1)^{z_1}e^{|u|-2|u|x}dx+$$
$$+ \int_{x\geq \frac{1}{|u|}+1} |u|^{k_2+z_2+1}|x|^{k_2}(x-1)^{z_2}e^{|u|-2|u|x}dx$$
$$\leq \int_1^\infty x^{k_1}(x-1)^{z_1}dx e^{-|u|x}|u|^{k_1+z_1+1} + \int_{\frac{1}{|u|}<x\leq \frac{1}{|u|}+1} |u|^{k_2+z_1+1}|x|^{k_2}(x-1)^{z_1}e^{|u|-2|u|x}dx+$$
$$+ \int_{x\geq \frac{1}{|u|}+1} |u|^{k_2+z_2+1}|x|^{k_2}(x-1)^{z_2}e^{|u|-2|u|x}dx.$$

The first integral converges under the condition

(69) $$z_1 + k_1 + 1 < 0,$$
(70) $$z_1 > -1.$$

Then

$$I_3 \leq B(-1-z_1-k_1, z_1+1)|u|^{k_1+z_1+1}e^{-|u|}+$$
$$+ (\gamma(z_1+1, 2) - \gamma(z_1+1, 2(1-|u|)))\frac{e^{-|u|}}{2^{z_1+1}}+$$
$$+ \Gamma(z_2+1, 2)\frac{(1+|u|)^{k_2}}{2^{z_2+1}}e^{|-u|}.$$

We can conclude that the some $I_1 + I_2 + I_3$ in this case of $|u| \leq 1$ can be bounded as

$$(g \bullet f)(u) \leq C|u|^{z_1+k_1+1}.$$

*2) Case $1 < |u| \leq 2$.*

$$I_1 = I_1^2 + I_1^4$$

$$= \int_{-\frac{1}{|u|} \leq x \leq 0} |u|^{k_1+z_2+1} |x|^{k_1} (1-x)^{z_2} e^{2|u|x-|u|} dx + \int_{x < -\frac{1}{|u|}} |u|^{k_2+z_2+1} |x|^{k_2} (1-x)^{z_2} e^{2|u|x-|u|} dx$$

$$= \int_0^2 x^{k_1} e^{-x} dx \frac{|u|^{z_2}}{2^{k_1+1}} e^{-|u|} + \int_2^\infty x^{k_2} e^{-x} dx \frac{(1+|u|)^{z_2}}{2^{k_2+1}} e^{-|u|}$$

$$= \gamma(k_1+1, 2) \frac{|u|^{z_2}}{2^{k_1+1}} e^{-|u|} + \Gamma(k_2+1, 2) \frac{(1+|u|)^{z_2}}{2^{k_2+1}} e^{-|u|},$$

$$I_2 = I_2^1 + I_2^2 + I_2^3$$

$$= \int_{1-\frac{1}{|u|} \leq x \leq \frac{1}{|u|}} |u|^{k_1+z_1+1} |x|^{k_1} (1-x)^{z_1} e^{-|u|} dx + \int_{0 < x < 1-\frac{1}{|u|}} |u|^{k_1+z_2+1} |x|^{k_1} (1-x)^{z_2} e^{-|u|} dx +$$

$$+ \int_{\frac{1}{|u|} < x < 1} |u|^{k_2+z_1+1} |x|^{k_2} (1-x)^{z_1} e^{-|u|} dx$$

$$\leq \int_{0 \leq x \leq 1} |u|^{k_1+z_1+1} |x|^{k_1} (1-x)^{z_1} e^{-|u|} dx + \int_{0 < x < 1-\frac{1}{|u|}} |u|^{k_1+1} |x|^{k_1} e^{-|u|} dx +$$

$$+ \int_{\frac{1}{|u|} < x < 1} |u|^{z_1+1} (1-x)^{z_1} e^{-|u|} dx$$

$$\leq B(k_1+1, z_1+1) |u|^{k_1+z_1+1} e^{-|u|} + \frac{(|u|-1)^{k_1+1}}{k_1+1} e^{-|u|} + \frac{(|u|-1)^{z_1+1}}{z_1+1} e^{-|u|},$$

and

$$I_3 = I_3^3 + I_3^4$$

$$= \int_{1 < |x| \leq 1+\frac{1}{|u|}} |u|^{k_2+z_1+1} |x|^{k_2} (x-1)^{z_1} e^{|u|-2|u|x} dx + \int_{x \geq 1+\frac{1}{|u|}} |u|^{k_2+z_2+1} |x|^{k_2} (x-1)^{z_2} e^{|u|-2|u|x} dx$$

$$\leq \int_{1 < |x| \leq 1+\frac{1}{|u|}} |u|^{k_2+z_1+1} (x-1)^{z_1} e^{|u|-2|u|x} dx + \int_{x \geq 1+\frac{1}{|u|}} |u|^{k_2+1} |x|^{k_2} e^{|u|-2|u|x} dx$$

$$\leq \gamma(z_1+1, 2) \frac{|u|^{k_2}}{2^{z_1+1}} e^{-|u|} + \Gamma(k_2+1, 2|u|+2) \frac{1}{2^{k_2+1}} e^{|u|}$$

Therefore, it is sufficient to bound the sum $I_1 + I_2 + I_3$ in this case as

$$(g \bullet f)(u) \leq C.$$

*3) Case $|u| > 2$.*

$$I_1 = I_1^2 + I_1^4$$
$$= \int\limits_{-\frac{1}{|u|} \leq x \leq 0} |u|^{k_1+z_2+1}|x|^{k_1}(1-x)^{z_2}e^{2|u|x-|u|}\mathrm{d}x + \int\limits_{x<-\frac{1}{|u|}} |u|^{k_2+z_2+1}|x|^{k_2}(1-x)^{z_2}e^{2|u|x-|u|}\mathrm{d}x$$
$$\leq \gamma(k_1+1, 2)\frac{|u|^{z_2}}{2^{k_1+1}}e^{-|u|} + \Gamma(k_2+1, 2)\frac{(|u|+1)^{z_2}}{2^{k_2+1}}e^{-|u|},$$
$$I_2 = I_2^2 + I_2^3 + I_2^4$$
$$= \int\limits_{0<x\leq\frac{1}{|u|}} |u|^{k_1+z_2+1}|x|^{k_1}(1-x)^{z_2}e^{-|u|}\mathrm{d}x + \int\limits_{1-\frac{1}{|u|}\leq x\leq 1} |u|^{k_2+z_1+1}|x|^{k_2}(1-x)^{z_1}e^{-|u|}\mathrm{d}x+$$
$$+ \int\limits_{\frac{1}{|u|}<x<1-\frac{1}{|u|}} |u|^{k_2+z_2+1}|x|^{k_2}(1-x)^{z_2}e^{-|u|}\mathrm{d}x$$
$$\leq \int\limits_{0<x\leq\frac{1}{|u|}} |u|^{k_1+1}|x|^{k_1}e^{-|u|}\mathrm{d}x + \int\limits_{1-\frac{1}{|u|}\leq x\leq 1} |u|^{z_1+1}(1-x)^{z_1}e^{-|u|}\mathrm{d}x+$$
$$+ \int\limits_{\frac{1}{|u|}<x\leq\frac{1}{2}} |u|^{k_2+z_2+1}|x|^{k_2}\frac{1}{2^{z_2}}e^{-|u|}\mathrm{d}x + \int\limits_{\frac{1}{2}<x<1-\frac{1}{|u|}} |u|^{k_2+z_2+1}\frac{1}{2^{k_2}}(1-x)^{z_2}e^{-|u|}\mathrm{d}x$$
$$\leq \frac{(|u|-1)^{z_2}}{k_1+1}e^{-|u|} + \frac{(|u|-1)^{k_2}}{z_1+1}e^{-|u|} - \frac{|u|^{z_2}}{2^{z_2}(k_2+1)}e^{-|u|} - \frac{|u|^{k_2}}{2^{k_2}(z_2+1)}e^{-|u|},$$
$$I_3 = I_3^3 + I_3^4$$
$$= \int\limits_{1<x\leq 1+\frac{1}{|u|}} |u|^{k_2+z_1+1}|x|^{k_2}(x-1)^{z_1}e^{|u|-2|u|x}\mathrm{d}x + \int\limits_{x>1+\frac{1}{|u|}} |u|^{k_2+z_2+1}|x|^{k_2}(x-1)^{z_2}e^{|u|-2|u|x}\mathrm{d}x$$
$$\leq \gamma(z_1+1, 2)\frac{|u|^{k_2}}{2^{z_1+1}}e^{-|u|} + \Gamma(z_2+1, 2)\frac{|u|^{k_2}}{2^{z_2+1}}e^{-|u|}.$$

Therefore, it is sufficient to bound the sum $I_1 + I_2 + I_3$ in this case as
$$(g \bullet f)(u) \leq C|u|^{\max\{k_2,z_2\}}e^{-|u|}.$$

We conclude that
$$(71) \qquad (g \bullet f)(u) \leq A_{k_1,k_2}^{z_1,z_2}|u|^{\zeta}e^{-|u|},$$
where $A_{k_1,k_2}^{z_1,z_2}$ is some constant depending on $k_1$, $k_2$, $z_1$ and $z_2$ and $\zeta = z_1+k_1+1$ for $|u| \leq 1$, $\zeta = \max\{z_2, k_2\}$ for $|u| > 1$.

Specifically, for $n \geq 2$,
$$(72) \qquad f^{\bullet n}(u) = \left(\prod_{i=1}^{n-1} A_{z_1,z_2}^{iz_1+i-1,z_2}\right)|u|^{\zeta}e^{-|u|},$$

where $\zeta = nz_1 + n - 1$ for $|u| \leq 1$, $\zeta = z_2$ for $|u| > 1$. The condition for convergence of integrals in (72) is

$$-(n-1) - nz_1 > 0 \implies z_1 < -\frac{n-1}{n}.$$

Finally,

$$(73) \qquad g \bullet f^{\bullet(2r)}(u) = \left( A_{k_1,k_2}^{(2r-1)z_1+2r-2, z_2} \prod_{i=1}^{2r-1} A_{z_1,z_2}^{(i-1)z_1+i-2, z_2} \right) |u|^\zeta e^{-|u|},$$

$\zeta = 2rz_1 + k_1 + 2r$ for $|u| \leq 1$, $\zeta = \max\{z_2, k_2\}$ for $|u| > 1$, the condition of convergence being

$$2r + 2rz_1 + k_1 < 0.$$

The claim follows.

## References

[1] J. B. Baillon, T. Cazenave, and M. Figuera. Essai sur le mouvement d'un fluide visqueux emplissant l'espace. *C. R. Acad. Sci., Paris*, 284:867–872, 1977.
[2] Yu. Bakhtin, E. I. Dinaburg, and Ya. G. Sinai. On solutions with infinite energy and enstrophy of the Navier-Stokes system. *Russian Math. Surveys*, 59(6):55–72, 2004.
[3] H. Berestycki and P. L. Lions. Existence d'ons solitaires dans des probléms nonlineares du type Klein-Gordon. *C. R. Paris*, 287:503–506, 1978.
[4] H. Berestycki, P. L. Lions, and L. A. Peletier. An ODE approach to the existence of positive solutions for semilinear problems in $\mathbb{R}^n$. *Ind. Uni. Math. J.*, 30:141–157, 1981.
[5] D. Gaidashev. Renormalization and a-priori bounds for Leray self-similar solutions to the generalized mild Navier-Stokes equations. *preprint: arXiv:2203.14646*, 2022.
[6] D. Gaidashev and A. Luque. Renormalization and existence of finite-time blow-up solutions for a one-dimensional analogue of the Navier-Stokes equations. *SIAM J. Math. Anal.*, 241(1):201–239, 2018.
[7] J. Ginibre and G. Velo. On a class of nonlinear Schrödinger equations: I. The Cauchy problem, general case. *J. Func. Analysis*, 32:1–32, 1979.
[8] R. T. Glassey. On the blowing up of solutions to the Cauchy problem for the nonlinear Schrödinger equation. *J. Math. Phys.*, 18:1794–1797, 1977.
[9] H. Hanche-Olsen, H. Holden, and E. Malinnikova. An improvement of the Kolmogorov-Riesz compactness theorem. *Expo. Math*, 37(1):84–91, 2019.
[10] Weinstein M. I. On the structure and formation of singularities in solutions to the nonlinear dispersive equations. *Comm. Partial Diff. Eq.*, 11:545–565, 1986.
[11] M. K. Kwong. Uniqueness of positive solutions of $\delta u - u + u^p = 0$ in $\mathbb{R}^n$. *Arch. Rat. Mech. An.*, 105:243–266, 1989.
[12] D. Li and Ya. G. Sinai. Blowups of complex solutions of the 3d-Navier-Stokes system and renormalization group method. *Journal of European Mathematical Society*, 10(2):267–313, 2008.
[13] D. Li and Ya. G. Sinai. Complex singularities of solutions of some 1D hydrodynamic models. *Physica D*, 237:1945–1950, 2008.
[14] D. Li and Ya. G. Sinai. Blowups of complex-valued solutions for some hydrodynamic models. *Regul. Chaotic Dyn.*, 15(4-5), 2010.
[15] D. Li and Ya. G. Sinai. Singularities of complex-valued solutions of the two-dimensional Burgers system. *J. Math. Phys.*, 51:521–531, 2010.

[16] F. Merle. On uniqueness and continuation properties after blow-up time of self-similar solutions on nonlinear Schrödinger equations with critical exponent and critical mass. *Commun. Pure Appl. Math.*, XLV:203–254, 1992.
[17] W. A. Strauss. Existence of solitary waves in higher dimensions. *Comm. Math. Phys.*, 55:149–162, 1977.
[18] W. A. Strauss. Everywhere defined wave operators. In *Nonlinear Evolution Equations*, pages 85–102. Academic Press, New York, 1978.
[19] W. A. Strauss. The nonlinear Schrödinger equation. In *Contemporary Developments in Continuum Mechanics and PDEs*. North-Holland, Amsterdam-New York - Oxford, 1978.
[20] Y. Tsutsumi. $L^2$-solutions for Nonlinear Schrödinger Equations and Nonlinear [g.

Department of Mathematics, Uppsala University, Uppsala, Sweden, gaidashev@math.uu.se